\documentclass[11pt]{article}
\usepackage{amsmath}
\usepackage{amsfonts}
\usepackage{amssymb}
\usepackage{graphicx}
\newcommand{\ls}[1]
   {\dimen0=\fontdimen6\the\font \lineskip=#1\dimen0
\advance\lineskip.5\fontdimen5\the\font \advance\lineskip-\dimen0
\lineskiplimit=.9\lineskip \baselineskip=\lineskip
\advance\baselineskip\dimen0 \normallineskip\lineskip
\normallineskiplimit\lineskiplimit \normalbaselineskip\baselineskip
\ignorespaces }

 \newcommand{\sfrac}[2]{\mbox{$\frac{#1}{#2}$}}

 \newcommand{\nfrac}[2]{\frac{\mbox{\normalsize{$#1$}}}{\mbox{\normalsize{$#2$}}}}
\def\squarebox#1{\hbox to #1{\hfill\vbox to #1{\vfill}}}

 \newcommand{\wtx}{\widetilde{x}}
 \newcommand{\psij}{\psi_{i,j}}
 \newcommand{\amat}{\left[\begin{array}{cccc}
                             a_1 &                                                  &        &                                  \\
                                 &                         a_2                      &        &  \hspace{-20pt}\mbox{\LARGE $0$} \\
                                 & \hspace{-15pt}\raisebox{-8pt}{\mbox{\LARGE $0$}} & \ddots &                                  \\
                                 &                                                  &        &              a_K                 \\
                                   \end{array} \right]\ \ }
\begin{document}
 \ls{1.5}
%
%
 \pagestyle{empty}
 \vspace*{2cm}

  \begin{center}
  \large
   THE OPTIMAL FILTERING OF MARKOV JUMP
  \medskip

   PROCESSES IN ADDITIVE WHITE NOISE
  \normalsize
  \bigskip
  \medskip

    by
  \bigskip

  M. Zakai
  \bigskip
  \bigskip
  \bigskip

  15 June 1965
  \bigskip

  Research Note No. 563
  \vspace{7cm}

  Applied Research Laboratory\\
  Sylvania Electronic System\\
  A Division of Sylvania Electric Products Inc.\\
  Waltham, Massachusetts 02154
 \end{center}
 \newpage
 \setcounter{page}{1}
%
%
 \pagestyle{plain}
  \begin{center}
   The optimal filtering of Markov jump

   processes in additive white noise.
  \end{center}
  \bigskip

  This note is based on Wonham \cite{Wonham}. The differences between this note and \cite{Wonham}
  are discussed in Section VIII.
  \medskip

  \noindent
  I.\ \ \underline{Statement of the Problem}.

   Let $x(t)$ be  Markov jump process with stationary transition probabilities and with a finite number of states.
   Let $a_1,a_2,\ldots,a_K$ be the states, let $p_{ij}(h)$ be the transition probabilities
   \[ p_{ij}(h) ={\rm Prob}\left\{x(t+h)=a_j\,\raisebox{-3pt}{\mbox{\huge $|$}}\,x(t)=a_i\right\}\,.\hspace{90pt} \]
   Let
   \[ p_{ij}(h)=\left\{  \begin{array}{cc}
                            1-\nu_ih+o(h) &\ \ \  i\!=\!j,\ \ h\,\mbox{\tiny $\searrow$}\,0  \\ \\
                            \nu_{ij}(h)+o(h) &\ \ \  i\!\ne\!j,\ \ h\,\mbox{\tiny $\searrow$}\,0 \\
                          \end{array}   \hspace{87pt}\right.\]

   \noindent
   where $\nu_i>0$,\ \ $\nu_{ij}\ge 0$ and $\displaystyle{\nu_i=\sum_{\stackrel{j=1}{i\ne j}}^K \nu_{ij}}\ \ \ ,i=1,\ldots,K$.
   \medskip

   Let $p_i(0)$ be the initial distribution of $x(0)$. In addition let $y(t)$ be a process given by
   \[ dy(t)=x(t)\,dt+\beta\,dw(t)\,,\ \ \ \ \ y(0)\!=\!0 \hspace{110pt} \]
   where $\beta$ is known and $w(t)$ is a standard Brownian motion. The problem is to find
  \[ p_j(t)={\rm Prob}\left\{x(t)\!=\!a_j\ \raisebox{-3pt}{\mbox{\huge $|$}}\ \pi_0^ty(\cdot)\right\}\hspace{125pt}\]
  where $\pi_a^by(\cdot)$ stands for $y(s),\,a\!\le\!s\!\le\!b$.
 \newpage
%
%
 \noindent
  II.\ \ \underline{A general expression for $p_j(t)$}

   The conditional probability\ \vspace{-4pt}\ ${\rm Prob}\left\{x(t)=a_j\ \raisebox{-3pt}{\mbox{\huge $|$}}\ A\right\}$\
   (where $A$ is some condition) is the same as the conditional expectation of the function $\delta_{ij}$\
   ($j$ fixed) given $A$. Applying Doob's theorem 8.8 (\cite[p.\,21]{Doob}),\ there exists a sequence $t_1,t_2,\ldots$,
   \ all in $[0,t]$,\ such that a.s.
   \[ p_j(t)={\rm Prob}\left\{ x(t)\!=\!a_j\ \raisebox{-3pt}{\mbox{\huge $|$}}\ y(t_1), y(t_2),\ldots\right\}\hspace{120pt}\]
   and by the martingale convergence theorem (\cite[Cor\,1 p.332]{Doob}) a.s.
   \begin{equation} \label{pjtlim}
    p_j(t)=\lim_{n\to\infty} \theta_j^{(n)}(t)\hspace{160pt}
   \end{equation}
   where\ \ \ \ \ \ \
   $\displaystyle{\theta_j^{(n)}(t)={\rm Prob}\left\{x(t)=a_j\ \raisebox{-3pt}{\mbox{\huge $|$}}\ y(t_1), y(t_2),\ldots,y(t_n)\right\}}$\,. \\
   In the following we will use:
   \[ {\rm Prob}\left\{x(t)\!\in\!A\,\raisebox{-3pt}{\mbox{\huge $|$}}\,\pi_0^ty(\cdot)\in B\right\}
          =\frac{{\rm Prob}\left\{x(t)\!\in\!A\,,\,\pi_0^ty(\cdot)\in B\right\}}
               {{\rm Prob}\left\{\pi_0^ty(\cdot)\in B\right\}}\ .   \hspace{30pt}              \]
   Consider now a fixed $t_\mu$,\ $0\!<\!t_\mu\!<\!t$;\ let $p_j(t\,|\,y(t_\mu))$ be the probability that $x(t)\!=\!a_j$ given $y(t_\mu)$
   and let $\xi_\mu$ be
   \[       \xi_\mu=\int_0^{t_\mu}\!x(t)\,dt\ .\hspace{100pt}\]
   Then
   \[ p_j\left(t\,\raisebox{-1pt}{\mbox{\Large$|$}}\,y(t_\mu)\right)\
         =\ \frac{\displaystyle{\sum_{i=1}^K p_i(0)\,p_{ij}(t)\int_{-\infty}^\infty\!\!\!P_1(y(t_\mu)\!-\!\xi_\mu)\,
                                                                                         P_2(\xi_\mu\,|\,x(0)\!=\!a_i,\,x(t)\!=\!a_j)\,d\xi_\mu}}
                           {\displaystyle{\sum_{j=1}^K (\mbox{The same expression as in the numerator})}}\ .\]
   Note that $P_1$ is normal $(0,\beta^2t_\mu)$. Let $\widetilde{x}(t)$ be a process independent of $x(t)$ and $w(t)$ with the same law as $x(t)$.
   Let $\widetilde{\xi}(s)=\int_{_0}^{^s}\!\!\widetilde{x}(t)\,dt$,\ \ then
   \newpage
%
%
  \[ p_j\left(t\,\raisebox{-1pt}{\mbox{\Large$|$}}\,y(t_\mu)\right)\
         =\ \frac{\displaystyle{\sum_{i=1}^K p_i(0)\,p_{ij}(t)\,
                    E\left\{\exp\left[-\frac{(y(t_\mu)\!-\!\widetilde{\xi}(t_\mu))^2}{2\beta^2t_\mu}\right]\,
                                          \ \raisebox{-3pt}{\mbox{\Huge $|$}}\ \widetilde{x}(0)\!=\!a_i,\,\widetilde{x}(t)\!=\!a_j\right\}}}
                 {\displaystyle{\sum_{j=1}^K\sum_{i=1}^K p_i(0)\,p_{ij}(t)\,
                    E\left\{\exp\left[-\frac{(y(t_\mu)\!-\!\widetilde{\xi}(t_\mu))^2}{2\beta^2t_\mu}\right]\,
                                          \ \raisebox{-3pt}{\mbox{\Huge $|$}}\ \widetilde{x}(0)\!=\!a_i,\,\widetilde{x}(t)\!=\!a_j\right\}}}\ \ . \]
   The conditioning of the expectations in the above expression are all $\widetilde{x}(s)$ paths which start at $s\!=\!0$ with $a_i$
   and terminate at $s\!=\!t$ in the state $a_j$. Now let $s_{r,n}=r\nfrac{t}{n},\ \ r=0,1,\ldots,n$.
   Let \vspace{-.5cm}
    \begin{eqnarray*}
     \eta_{r,n}&=&y(s_{r+1,n})-y(s_{r,n})\hspace{160pt}\\
     \widetilde{\xi}_{r,n}&=&\int_{_{r\frac{t}{n}}}^{^{(r+1)\frac{t}{n}}}\hspace{-17pt}\widetilde{x}(s)\,ds\ \ .
    \end{eqnarray*}
   Then by the same arguments as above
   \begin{eqnarray}
    \lefteqn{p_j\left(\!t\,\raisebox{-1pt}{\mbox{\Large $|$}}\,y(r\!\sfrac{t}{n}),\,r\!=\!0,1,\ldots,n\!\right)=} \nonumber\\
        &&=\ \frac{\displaystyle{\sum_{i=1}^K p_i(0)\,p_{ij}(t)\,
                    E\!\left\{\!\exp\!-\!\!\sum_{r=0}^{n-1}\frac{(\eta_{r,n}\!-\!\widetilde{\xi}_{r,n})^2}{2\beta^2t/n}\,
                                          \raisebox{-3pt}{\mbox{\Huge $|$}}\,\widetilde{x}(0)\!=\!a_i,\,\widetilde{x}(t)\!=\!a_j\!\right\}}}
                 {\displaystyle{\sum_{j=1}^K (\mbox{numerator})}}\ \ .   \label{pjratio}
   \end{eqnarray}
   The argument of the exponential is
     \[ \sum_{r=0}^{n-1} \left(\eta_{r,n}^2-2\eta_{r,n}\widetilde{\xi}_{r,n}+\widetilde{\xi}_{r,n}^{\,2}\right)\,\frac{1}{2\beta^2t/n}\ .\]
   The first term will be cancelled by the same term in the denominator. The last term converges a.s. as $n\to\infty$ to
   \[  \frac{1}{2\beta^2}\,\int_0^t \widetilde{x}(t)^2\,dt\ \hspace{160pt}.\]
   The middle term converges a.s. to
   \[  -\frac{1}{\beta^2}\,\int_0^t\!\widetilde{x}(t)\,dy(t).   \hspace{165pt} \]
 \newpage
%
%
   \noindent
   We want to apply these results to the evaluation of the limit of the numerator of (\ref{pjratio}) as $n\!\to\!\infty$.
   In order to do that we have to show that if $f_n\!\to\!f$ as $n\!\to\!\infty$, then $E(f_n\,|\ \,)\to E(f\,|\ \,)$. Since
   \begin{eqnarray*}
    \lefteqn{\exp -\left[\sum_{r=0}^{n-1}\left(-2\eta_{r,n}\,\widetilde{\xi}_{r,n}+\widetilde{\xi}_{r,n}^{\ 2}\right)\right]}\\
        &&\hspace*{50pt}\le \exp\left[\,A^2_{\rm max}\!\cdot\!t+2A_{\rm max}\!\cdot\!\sup\{|y(t_1)\!-\!y(t_2)|,\,t_1,t_2\!\in\![0,t]\}\,\right]
   \end{eqnarray*}
   it follows, by dominated convergence that the limit of the numerator of (\ref{pjratio}), as $n\to\infty$, is
   \begin{equation} \label{psij}
    \psi_j(t)=\sum_{i=1}^K p_i(0)\,p_{ij}(t)\,E\left\{\exp\left[-\frac{1}{2\beta^2}\!\int_0^t\!\!\widetilde{x}^{\,2}(s)\,ds
                                                                   +\frac{1}{\beta^2}\!\int_0^t\!\!\widetilde{x}(s)\,dy(s)\right]\,
                                  \raisebox{-3pt}{\mbox{\Huge $|$}}\, \widetilde{x}(0)\!=\!a_i,\,\widetilde{x}(t)\!=\!a_j\right\}
   \end{equation}
   and the conditioning is with respect to all the paths which start at $\widetilde{x}(0)\!=\!a_i$\ and terminate at $\widetilde{x}(t)\!=\!a_j$.
   Similarly the limit of the denominator is $\sum_{j=1}^K\!\psi_j(t)$\ where $\psi_j(t)$ is given by equation (\ref{psij}).
   Since $\psi_j(t)\!>\!0$ a.s. we have
   \begin{equation} \label{limpj}
    \lim_{n\to\infty} p_j\left(t\ |\ y(r\sfrac{t}{n}\right),\ r\!=\!0,1,\ldots,n)=\frac{\psi_j(t)}{\sum_{i=1}^K \psi_i(t)}\ .\hspace{30pt}
   \end{equation}

   The limits (\ref{psij}) and (\ref{limpj}) were obtained by a particular sequence of partitions of $[0,t]$, but it is clear
   that the same result will hold for any sequence of partitions $\{s_{r,n}\}$ such that $0\!=\!s_{0,n}\!<\!s_{1,n}\!<\!\ldots\!<\!s_{n,n}\!=\!t$
   and such that $\max\hspace{-13pt}\raisebox{-8pt}{$r$}\hspace{13pt}(s_{r+1,n}\!-\!s_{r,n})\!\to\!0$ as $n\!\to\!\infty$. We may therefore
   use a sequence for which (\ref{pjtlim}) is true. Therefore
   \begin{equation} \label{pjtratio}
    p_j(t)=\frac{\psi_j(t)}{\sum_{i=1}^K\psi_i(t)}\ . \hspace{160pt}
   \end{equation}
 \newpage
%
%
 \noindent
   III.\ \ \underline{The stochastic differential equation for $\psi_j(t)$}.
   \bigskip

   Let $\psi_{i,j}(a,b),\ \ (b\!>\!a)$ be
   \[ \psi_{i,j}(a,b)=p_{ij}(b\!-\!a)\,E\left\{\exp\left[-\frac{1}{2\beta^2}\!\int_a^b\!\!\widetilde{x}^{\,2}(s)\,ds
                                                                   +\frac{1}{\beta^2}\!\int_a^b\!\!\widetilde{x}(s)\,dy(s)\right]\,
                                  \raisebox{-3pt}{\mbox{\Huge $|$}}\, \widetilde{x}(a)\!=\!a_i,\,\widetilde{x}(b)\!=\!a_j\right\}\ . \]
 Then, comparing with~(\ref{psij}):
 \[ \psi_j(t)=\sum_{i=1}^K p_i(0)\,\psi_{i,j}(0,t)\ .\hspace{135pt}\]
 Consider a fixed realization of $\pi_0^{t+h}y(\cdot)$, we prove now that
 \begin{equation} \label{psijtplush}
  \psi_j(t\!+\!h)=\sum_{i=1}^K\psi_i(t)\cdot\psi_{i,j}(t,t\!+\!h)\ .\hspace{100pt}
 \end{equation}

  \noindent
  Proof:
 \begin{eqnarray*}
  \lefteqn{\psi_{i,j}(0,t\!+\!h)}\\
            &&=p_{ij}(t\!+\!h)\,E\left\{\exp\left[-\!\int_0^{t+h}\!\!\!\!\!\!\!\ldots\ +\!\!\int_0^{t+h}\!\!\!\!\!\!\ldots\ \ \right]\,
                                \raisebox{-3pt}{\mbox{\Huge $|$}}\,\widetilde{x}(t\!+\!h)\!=\!a_j,\,\widetilde{x}(0)\!=\!a_i\right\} \\
            &&=p_{ij}(t\!+\!h)\,\sum_{k=1}^KE\left\{\exp\left[-\!\int_0^{t+h}\!\!\!\!\!\!\!\ldots\ +\!\!\int_0^{t+h}\!\!\!\!\!\!\ldots\ \right]\,
      \raisebox{-3pt}{\mbox{\Huge $|$}}\,\widetilde{x}(t\!+\!h)\!=\!a_j,\,\widetilde{x}(t)\!=\!a_k,\,\widetilde{x}(0)\!=\!a_i\right\} \\
     &&\hspace*{175pt}\mbox{\LARGE $\cdot$}\ {\rm Prob}\left\{\widetilde{x}(t)\!=\!a_k\ \mbox{\Large $|$}
                                           \ \widetilde{x}(t\!+\!h)\!=\!a_j,\,\widetilde{x}(0)\!=\!a_i\right\}\ .
 \end{eqnarray*}
 Since $\widetilde{x}(t)$ is a Markov process, the conditional expectation becomes the product of two conditional expectations
 (since, given $\widetilde{x}(t),\ \ \widetilde{x}(t\!-\!\alpha)$ and $\widetilde{x}(t\!+\!\beta)$ are independent for $\alpha,\,\beta\!>\!0$).
 Moreover
 \[ {\rm Prob}\,\left\{\widetilde{x}(t)=a_k\ \mbox{\Large $|$}\ \widetilde{x}(t\!+\!h)\!=\!a_j,\,\widetilde{x}(0)\!=\!a_i\right\}
                       =\frac{p_{ik}(t)\,p_{kj}(h)}{p_{ij}(t\!+\!h)}\ . \]
 \newpage
%
%
 Therefore
 \begin{eqnarray*}
  \psij(0,t\!+\!h)&=&\sum_{k=1}^Kp_{ik}(t)\,E\left\{\exp\left[-\!\int_0^t\!\!\!\ldots\ +\!\!\int_0^t\!\!\!\ldots\ \
                                    \,\raisebox{-3pt}{\mbox{\Huge $|$}}\,\wtx(t)\!=\!a_k,\ \wtx(0)\!=\!a_i\right] \right\}\\
         &&\hspace{20pt}\mbox{\Large $\cdot$}\ p_{kj}(h)\,E\left\{\exp\left[-\!\int_0^t\!\!\!\ldots\ +\!\!\int_0^t\!\!\!\ldots\ \
                                                     \raisebox{-3pt}{\mbox{\Huge $|$}}\,\wtx(t\!+\!h)\!=\!a_j,\ \wtx(t)\!=\!a_k\right] \right\}\\
         &=&\sum_{k=1}^K\psi_{i,k}(0,t)\,\mbox{\Large $\cdot$}\,\psi_{k,j}(t,t\!+\!h)
 \end{eqnarray*}
 which is the required result.
 \bigskip

 \noindent
 Since the $x(t)$ process is Markov and the $w(t)$ process has independent increments, it follows that $\psi_{i,k}(0,t)$
 and $\psi_{k,j}(t,t\!+\!h)$ are conditionally independent given $x(t)$.
 \underline{Therefore the process} \underline{$\left(x(t),\psi_1(t),\psi_2(t),\ldots,\psi_K(t)\right)$ is a $K\!+\!1$ dimensional Markov process}.
 \medskip

 $\psij(t,t\!+\!h)$ will now be evaluated for small $h$. Assuming that $h$ is small enough so that the possibility that more than one transition
 in $[t,t\!+\!h]$ can be ignored we have
 \begin{eqnarray*}
  \psi_{i,i}(t,t\!+\!h)&\cong&(1\!-\!\nu_i\,h)\,\exp\left[-\frac{a_i^2}{2\beta^2}\,h+\frac{a_i}{\beta^2}\,(y(t\!+\!h)\!-\!y(t))\right]\\
  \psi_{i,j}(t,t\!+\!h)&\cong&\nu_{ij}\,h\,\exp\left[-\frac{a_j^2\,h\,\theta_1}{2\beta^2}+\frac{a_j\,\theta_2}{\beta^2}\,(y(t\!+\!h)\!-\!y(t))\right]
                                                                                                                    \hspace{30pt}j\!\ne\!i\\
 \end{eqnarray*}
 The factors $\theta_1$ and $\theta_2$ were included in the last expression in order to indicate that it is unknown where in $[t,t\!+\!h]$ the
 transition occurred; it will turn out that this is immaterial. Setting now $y(t\!+\!h)\!-\!y(t)=\int_t^{t+h}\!\!x(s)\,ds+\beta\,w(t\!+\!h)-\beta\,w(t)$
 and expanding the exponential in a power series we obtain
 \[ \psi_{i,i}(t,t\!+\!h)\cong1-\nu_ih-\frac{a_i^2}{2\beta^2}\,h+\frac{a_i}{\beta^2}\,(y(t\!+\!h)\!-\!y(t))+\frac{a_i^2}{2\beta^4}\,\beta^2
      \,(w(t\!+\!h)\!-\!w(t))^2+o_1(h,(\Delta_hw)^2)  \]
 where $o_1(h,(\Delta_hw)^2)$ denotes the terms omitted. Also
 \newpage
%
%
 \[  \psij(t,t\!+\!h)\cong\nu_{ij}\,h+o_2(h,(\Delta_hw)^2)\ .\hspace{175pt}\]
 Substituting into~(\ref{psijtplush}) we have
 \begin{eqnarray*}
   \psi_j(t\!+\!h)&=&\psi_j(t)-\nu_j\,\psi_j(t)\,h-\frac{a_j^2}{2\beta^2}\,h\,\psi_j(t)+\frac{a_j\,\psi_j(t)}{\beta^2}\,(y(t\!+\!h)\!-\!y(t))\\
         &&\hspace{10pt}+\,\frac{a_j^2}{2\beta^2}\,(w(t\!+\!h)\!-\!w(t))^2+\sum_{\stackrel{i=1}{i\ne j}}^K\psi_i(t)\,\nu_{ij}\,h+o_3(h,(\Delta_hw)^2)\ .
 \end{eqnarray*}
 \vspace{-.5cm}

 \noindent
 Therefore
 \vspace{-.2cm}
 \begin{eqnarray}
  \psi_j(t)-\psi_j(0)
      &\!\!=\!\!&\mbox{\Huge $\int$}_{\!\!\!\!0}^t\left[-\nu_j\,\psi_j(s)
                                      +\sum_{\stackrel{i=1}{i\ne j}}^K\psi_i(s)\,\nu_{ij}-\frac{a_j^2\,\psi_j(s)}{2\beta^2}\right]\,ds \nonumber\\
      &&\hspace*{10pt}+\,\int_0^t\frac{a_j\,\psi_j(s)}{\beta^2}\,dy(s)
                       +\lim_{h\to 0}\ \sum_{r=1}^{t/h} \frac{a_j^2\,\psi_j(hr)}{2\beta^2}\,\left(\,w((r\!+\!1)h)\!-\!w(rh)\,\right)^2 \nonumber\\
      &&\hspace*{10pt}+\,\lim_{h\to 0}\ \sum_{r=1}^{t/h} o_3(h,(\Delta_hw)^2)\ . \label{psijlimit}
 \end{eqnarray}
  The first sum can be shown~\cite{WZ2} to converge a.s. to
  \[ \int_0^t \frac{a_j^2\,\psi_j(s)}{2}\,ds\ .\hspace{150pt}\]
  The second sum can be shown~\cite{Ito1} to converge to $0$. (For example, if $|f(s)|\!\le\!M$ in $[0,t]$ then the sum
  \ $\sum_{r=1}^{t/h}\!\!f(rh)\,h\,(w((r\!+\!1)h)\!-\!w(rh))$ is bounded by
  \[ Mt\,\sup_t\,|w(t\!+\!h)\!-\!w(t)| \]
  but since $w(t)$ is a.s. continuous it is also uniformly continuous and the last term converges to zero as $h\!\to\!0$,
  hence the sum converges a.s. to zero).\vspace{.1cm}

  \noindent
  Equation~(\ref{psijlimit}) becomes:\vspace{-.5cm}
  \begin{equation} \label{dpsij}
   d\psi_j(t)=-\nu_j\,\psi_j(t)\,dt+\sum_{\stackrel{i=1}{i\ne j}}^K \nu_{ij}\,\psi_i(t)\,dt+\frac{a_j\,\psi_j(t)}{\beta^2}\,dy(t),
  \end{equation}
 \newpage
%
%
 \noindent with the initial condition $\psi_j(0)\!=\!p_j(0)$.
 \medskip

 \noindent
 Let $\phi(t)=\displaystyle{\sum_{j=1}^K \psi_j(t)}$.\ Then
  \[ d\phi(t)=\sum_{i=1}^K \frac{a_i\,\psi_i(t)}{\beta^2}\,dy(t)\,,\hspace{180pt}\]
  and \vspace{-.1cm}
  \begin{equation} \label{pjt}
   p_j(t)=\frac{\psi_j(t)}{\sum_{i=1}^K\psi_i(t)}=\frac{\psi_j(t)}{\phi(t)}\ . \hspace{165pt}
  \end{equation}

  \noindent
  Equation~(\ref{dpsij}) is a stochastic differential equation for $\psi_j(t)$ from which $p_j(t)$ can be obtained by~(\ref{pjt}).

  The Langevin equation\footnote{see~\cite{Ito1} or~\cite{WZ3}} corresponding to~(\ref{dpsij}) can be derived using equation (4.30) of~\cite{WZ1};
  the result is
  \begin{eqnarray}
   \frac{d\psi_j{t}}{dt}&=&-\nu_j\,\psi_j(t)+\sum_{i\ne j} \nu_{ij}\,\psi_i(t)+\frac{1}{2}\,\frac{a_j^2\,\psi_j(t)}{\beta^2}
                   +\frac{a_j\,\psi_j(t)}{\beta^2}\,\frac{dy(t)}{dt} \label{langevin}
  \end{eqnarray}
  and $dy(t)/dt$ is $x(t)$ plus ``white noise".
 \bigskip

 \noindent
   IV.\ \ \underline{The stochastic differential equation for $p_j(t)$}.

 Since (by definition of $\psi_j$ and $\phi$)\ \ \ $\phi(t)\!\ne\!0$ a.s.,\ we may apply Ito's rule of differentiation~\cite{Ito1} to~(\ref{pjt}):
 \begin{eqnarray*}
  dp_j(t)&=&\frac{d\psi_j(t)}{\phi(t)}-\frac{\psi_j(t)\,d\phi(t)}{\phi^2(t)}
                         -\frac{1}{\phi^2(t)}\,\frac{a_j\,\psi_j(t)}{\beta^2}\,\left(\sum_{i=1}^K\frac{a_i\psi_i(t)}{\beta^2}\right)\beta^2dt\\
          &&\hspace*{20pt}+\ \frac{\psi_j(t)}{\phi^3(t)}\,\left(\sum_{i=1}^K\frac{a_i\psi_i(t)}{\beta^2}\right)^{\!2}\!\beta^2dt\ .
 \end{eqnarray*}
 Substituting for $d\psi_j$ and $d\phi$, and setting $\frac{\sum_{j=1}^K a_j\,\psi_j(t)}{\phi(t)}=\overline{x}(t)$ we get
 \newpage
%
%
  \begin{eqnarray}
   dp_j(t)&=&-\nu_j\,p_j(t)\,dt+\sum_{i\ne j}^Kp_i(t)\,\nu_{ij},dt \nonumber\\
           &&\hspace*{20pt}+\ \beta^{-2}\left(a_j\!-\!\overline{x}(t)\right)\,p_j(t)\,dy(t) \label{dpjt}\\
           &&\hspace*{20pt}+\ \beta^{-2}\,\overline{x}(t)\,p_j(t)\left(a_j\!-\!\overline{x}(t)\right)\,dt
                 \hspace{35pt}[\,\mbox{\small $\overline{x}(t)\!=\!\sum_{j=1}^Ka_jp_j(t)$}\,]\nonumber
  \end{eqnarray}
 which are the equations derived by Wonham. Note that as $\beta^2\!\to\!\infty$ we get the Kolmogorov forward equation (as expected).
 The equations for $\psi_j(t)$ are more elegant than those for $p_j(t)$. However, the equations for $p_j(t)$ are probably more useful
 for applications since they ensure that the output is always in $[0,1]$ (while $\psi_j(t)$ can be anywhere in $(0,\infty)$).
 It seems also that perhaps the $p_j(t)$ may have a stationary distribution while $\psi_j(t)$ may not have such a distribution.

 The Langevin equation corresponding to~(\ref{dpjt}) is, by eq.~(4.30) of~\cite{WZ1},
  \begin{eqnarray}
   \frac{dp_j(t)}{dt}&=&-\nu_j\,p_j+\sum_{i\ne j}p_i\,\nu_{ij}+\frac{1}{2}\,p_j\,\beta^{-2}\,\left(a_j^2-\sum_{i=1}^Ka_i^2p_i\right)\hspace{80pt}\nonumber\\
                      &&\hspace*{20pt}+\ \beta^{-2}\left(a_j-\overline{x}\right)\,p_j\,\frac{dy}{dt}\ .
  \end{eqnarray}

 \noindent
 V.\ \ \underline{Example - The random telegraph signal}.
   \smallskip

  In this case
   \vspace{-.83cm}
  \begin{eqnarray*}
   a_1&\!\!=\!\!&1,\ \ \ \ a_2=-1\hspace{150pt}\\
   \nu_i&\!\!=\!\!&\nu_{ij}=\nu\,;\ i,j\!=\!1,2
  \end{eqnarray*}
  where $\nu$ is the expected number of jumps.

  \noindent
  Let  \hspace{27pt} $q(t)=p_1(t)-p_2(t)$\\
  then \hspace{22pt} $\overline{x}(t)\!=\!q(t)$\ \ \ \ \ \  and the equations of the last paragraph become
  \[ \hspace*{25pt}dq(t)=-2\nu\,q(t)\,dt-\beta^{-2}q(t)(1-q^2\!(t))\,dt+\beta^{-2}(1\!-\!q^2\!(t))\,dy(t)\hspace{55pt} \]
  or, equivalently
  \begin{eqnarray*}
   dq(t)&\!\!=\!\!&-2\nu\,q(t)\,dt-\beta^{-2}q(t)(1-q^2\!(t))\,dt+\beta^{-2}(1\!-\!q^2\!(t))\,x(t)\,dt\hspace{18pt}\\
       &&\hspace*{15pt}+\ \beta^{-1}(1\!-\!q^2\!(t))\,dw(t)\ .
  \end{eqnarray*}
 \newpage
%
%
 \noindent
 The Langevin equivalent of this stochastic differential equation is given by (\cite{WZ1},~\cite{WZ3})
  \[ \frac{dq(t)}{dt}=-2\nu\,q(t)+\beta^{-2}(1\!-\!q^2\!(t))\,x(t)+\beta^{-2}(1\!-\!q^2\!(t))\,n(t) \]
 where $n(t)$ is ``white noise". Let $x(t)\!+\!n(t)\!=\!r(t)$, then we have the Riccati equation:
  \[ \frac{dq(t)}{dt}=-2\nu\,q(t)+\beta^{-2}(1\!-\!q^2\!(t))\,r(t) \hspace{35pt}[\,\mbox{\small $r(t)\!=\!x(t)\!+\!n(t)$}\,]\ .\]
  The physical filter to compute $q(t)$ will therefore be

 
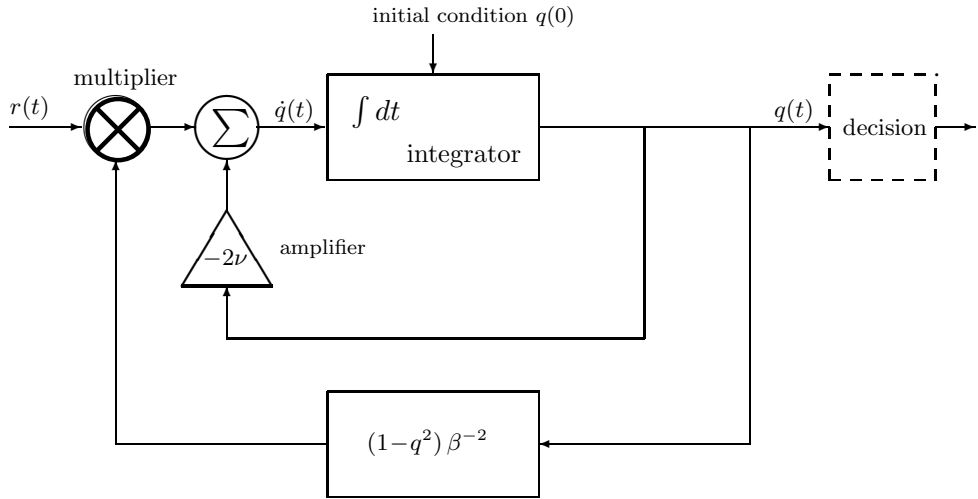
\begin{picture}(200,200)(0,0)
  \put(30,20){\line(1,0){80}}
  \put(110,0){\line(0,1){40}}
  \put(110,0){\line(1,0){80}}
  \put(110,40){\line(1,0){80}}
  \put(190,0){\line(0,1){40}}
      \put(125,18){\mbox{\footnotesize $(1\!-\!q^2)\,\beta^{-2}$}}
  \put(270,20){\vector(-1,0){80}}
  \put(30,20){\vector(0,1){108}}
  \put(-10,140){\vector(1,0){27}}
     \put(-10,145){\mbox{\footnotesize $r(t)$}}
  \put(30,140){\circle{25}}
      \put(17,133){\mbox{\Huge  $\bigotimes$}}
      \put(14,156){\mbox{\footnotesize \rm multiplier}}
  \put(42,140){\vector(1,0){18}}
  \thicklines
  \put(72,140){\circle{25}}
      \put(66,136){\mbox{\Large  $\mathbf \sum$}}
  \thinlines
  \put(83,140){\vector(1,0){26}}
     \put(90,144){\mbox{\footnotesize $\dot{q}(t)$}}
  \put(110,120){\line(0,1){40}}
  \put(110,120){\line(1,0){80}}
  \put(110,160){\line(1,0){80}}
  \put(190,120){\line(0,1){40}}
     \put(119,143){$\int dt$}
     \put(140,127){\small \rm integrator}
  \put(150,175){\vector(0,-1){15}}
       \put(127,180){\mbox{\scriptsize initial condition $q(0)$}}
  \put(190,140){\vector(1,0){110}}
     \put(279,144){\mbox{\footnotesize $q(t)$}}
  \multiput(300,120)(8,0){5}{\line(1,0){4}}
  \multiput(300,160)(8,0){5}{\line(1,0){4}}
  \multiput(300,120)(0,8){5}{\line(0,1){4}}
  \multiput(340,120)(0,8){5}{\line(0,1){4}}
      \put(305,137){\footnotesize \rm decision}
  \put(339,157){$\cdot$}
  \put(340,140){\vector(1,0){15}}
   \put(270,20){\line(0,1){120}}
  \put(230,140){\line(0,-1){80}}
  \put(230,60){\line(-1,0){158}}
  \put(72,60){\vector(0,1){20}}
    \thicklines
  \put(55,80){\line(1,0){34}}
  \put(55,80){\line(3,5){17}}
  \put(89,80){\line(-3,5){17}}
      \put(63,88){\mbox{\footnotesize $-2\nu$}}
      \put(92,92){\scriptsize \rm amplifier}
    \thinlines
  \put(72,108){\vector(0,1){20}}
 \end{picture}
 \medskip\medskip

 \noindent or:
  
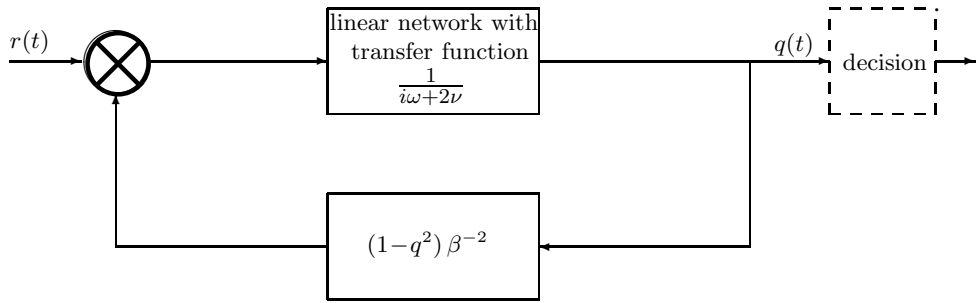
\vspace{-.5cm}

 \begin{picture}(200,200)(0,0)
  \put(30,70){\line(1,0){80}}
  \put(110,50){\line(0,1){40}}
  \put(110,50){\line(1,0){80}}
  \put(110,90){\line(1,0){80}}
  \put(190,50){\line(0,1){40}}
      \put(125,68){\mbox{\footnotesize $(1\!-\!q^2)\,\beta^{-2}$}}
  \put(270,70){\vector(-1,0){80}}
  \put(30,70){\vector(0,1){58}}
  \put(-10,140){\vector(1,0){27}}
     \put(-10,145){\mbox{\footnotesize $r(t)$}}
  \put(30,140){\circle{25}}
      \put(17,133){\mbox{\Huge  $\bigotimes$}}
  \put(42,140){\vector(1,0){68}}
  \put(110,120){\line(0,1){40}}
  \put(110,120){\line(1,0){80}}
  \put(110,160){\line(1,0){80}}
  \put(190,120){\line(0,1){40}}
      \put(111,152){\footnotesize \rm linear network with}
      \put(119,141){\footnotesize \rm transfer function}
      \put(136,128){\mbox{\normalsize $\frac{1}{i\omega+2\nu}$}}
  \put(190,140){\vector(1,0){110}}
     \put(279,144){\mbox{\footnotesize $q(t)$}}
  \multiput(300,120)(8,0){5}{\line(1,0){4}}
  \multiput(300,160)(8,0){5}{\line(1,0){4}}
  \multiput(300,120)(0,8){5}{\line(0,1){4}}
  \multiput(340,120)(0,8){5}{\line(0,1){4}}
      \put(305,137){\footnotesize \rm decision}
  \put(339,157){$\cdot$}
  \put(340,140){\vector(1,0){15}}
   \put(270,70){\line(0,1){70}}
%
 \end{picture}
  \newpage
%
%
 If, instead of the analog filter we use a digital computer we have to distinguish between two cases. Let $\omega_n$ be the cutoff frequency
 of the ``white noise" and $\omega_s$ be the sampling frequency of the computer. Case 1: $\nu\!\ll\!\omega_n\!\ll\!\omega_s$,\ Case 2:
 $\nu\!\ll\!\omega_s\!\ll\!\omega_n$\,. It follows from~\cite{WZ1} and~\cite{WZ3} that in case 1 the computer should be programmed to solve
 the Langevin equation. In case 2, Maruyama's approximation theorem is applicable~\cite{WZ1} and the computer should be programmed to solve
 Ito's equation (via Maruyama's approximation).

 An error analysis for this example is discussed in~\cite{Wonham}.
 \bigskip

 \noindent
 VI.\ \ \underline{Some transformations on $\psi_j(t)$}.
   \smallskip

  Equation~(\ref{langevin}) can be rewritten as
  \[ \hspace*{10pt}\frac{d\psi(t)}{dt}=A\,\cdot\,\psi(t)+\amat\!\!\frac{(x(t)\!+\!n(t))}{\beta^2}\,\cdot\,\psi(t)
                                                                                        \hspace{120pt}(10{\rm a}) \]
  where $A$ is a constant matrix and $\psi$ is the vector $(\psi_1,\ldots,\psi_K)^T$. \\Setting
  \[\Gamma(t)=e^{-At}\cdot\psi(t)\hspace{210pt}\]
  we get the Langevin equation for $\Gamma(t)$:
  \begin{equation} \label{langevingamma}
   \frac{d\Gamma(t)}{dt}=e^{-At}\amat\!\!e^{At}\ \ \frac{(x(t)\!+\!n(t))}{\beta^2}\ \Gamma(t)\ .\hspace{115pt}
  \end{equation}

  Setting $\theta_j(t)\!=\!\log\psi_j(t)$\\
  hence
  \begin{equation} \label{theta}
   \frac{d\theta_j(t)}{dt}=-\nu_j+\frac{1}{2}\,\frac{a_j^2}{\beta^2}\,+\,\sum_{i\ne j}\nu_{ij}\,e^{\theta_i(t)\!-\!\theta_j(t)}
                                                             +\frac{a_j}{\beta^2}(x(t)\!+\!n(t))\hspace{32pt}
  \end{equation}
  and since $\theta_i$ does not appear in front of the last term, this is the Ito as well as the Langevin equation for $\theta_j(t)$\,.
  \newpage
%
%
  \noindent VII.\ \ \underline{The prediction of $x(t)$}.

   \noindent
   The problem is now to find the probability that $x(t\!+\!h)\!=\!a_j,\ h\!>\!0$,\ conditioned on $\pi_0^ty(\cdot)$\,. The result follows directly from
  \begin{eqnarray*}
   \lefteqn{{\rm Prob}\left\{x(t\!+\!h)\!\in\!A\,\raisebox{-3pt}{\mbox{\Huge $|$}}\,\pi_0^ty(\cdot)\!\in\!B\right\}=}\\
     &&=\sum_{i=1}^K{\rm Prob}\left\{x(t\!)\!=\!a_i\,\raisebox{-3pt}{\mbox{\Huge $|$}}\,\pi_0^ty(\cdot)\!\in\!B\right\}\ \cdot\
                    {\rm Prob}\left\{x(t\!+\!h)\!\in\!A\,\raisebox{-3pt}{\mbox{\Huge $|$}}\,x(t)\!=\!a_i,\ \pi_0^ty(\cdot)\!\in\!B\right\}\ .
  \end{eqnarray*}
  Since\ ${\rm Prob}\left\{x(t\!+\!h)\,\raisebox{-3pt}{\mbox{\Huge $|$}}\,x(t),\ \pi_0^ty(\cdot)\right\}
          ={\rm Prob}\left\{x(t\!+\!h)\,\raisebox{-3pt}{\mbox{\Huge $|$}}\,x(t)\right\}$\,,\ it follows that (as expected)
  \begin{equation} \label{pred}
   {\rm Prob}\left\{x(t\!+\!h)\!=\!a_j\,\raisebox{-3pt}{\mbox{\Huge $|$}}\,\pi_0^ty(\cdot)\right\}=\sum_{i=1}^Kp_i(t)\,p_{ij}(h)\hspace{60pt}
  \end{equation}
  where $p_i(t)$ is the solution to~(\ref{dpjt}) and $p_{ij}(h)$ are defined at the beginning of section I. The extension of~(\ref{pred}) to
  the probability distribution of functionals on $x(s),\ s\!\ge\!t$,\ conditioned on $\pi_0^ty(\cdot)$ is obvious.
  \bigskip

  \noindent VIII.\ \ \underline{Remarks}.

   Section~I,~II and the first halves of~III and~V follow from Wonham~\cite{Wonham}. Instead of proceeding directly to obtain the
   stochastic differential equations for $p_j(t)$ as done in~\cite{Wonham} we first derive the stochastic differential equations for $\psi_j(t)$
    (section~III) from which the stochastic differential equations for $p_j(t)$ are derived by a singular transformation (section~IV).
    The equations for $\psi_j(t)$ are considerably simpler and are of a standard form (section~VI). It is believed that a similar approach
    can be used in the case treated by Stratonovich and Kushner (where $x(t)$ is a diffusion process). The treatment in this note is restricted
    to $\beta\!=$\,const, the treatment in~\cite{Wonham} is for $\beta\!=\!\beta(t)$ where $\beta(t)$ is continuously differentiable and bounded
    away from zero. The extension of the arguments and results of this note to $\beta\!=\!\beta(t)$ is straightforward.
 \newpage
%
%
 A question which was left open in~\cite{Wonham} was the problem of the realization of the results as physical ``filters".
 Recent work reated to this problem~\cite{WZ1},\ ~\cite{WZ3} gives answers to this question. Equation~(4.30) of~\cite{WZ1}, which was used in this note,
 was derived in~\cite{WZ1} by a heuristic argument. Unpublished calculations (for piecewise linear approximations to the Brownian motion) show
 that~(4.30) is correct. A short discussion on the realization problem is included in section~V.
 \bigskip\bigskip

\end{document}